\documentclass[11pt]{amsart}
\usepackage{amssymb}
\usepackage{epsfig}
\usepackage{graphicx}
\newtheorem{theorem}{Theorem}
\newtheorem{lemma}[theorem]{Lemma}

\newtheorem{conjecture}[theorem]{Conjecture}

\newtheorem{corollary}[theorem]{Corollary}
\newtheorem{problem}[theorem]{Problem}
\newtheorem{definition}[theorem]{Definition}

\newcommand{\Z}{\mathbb{Z}}

\newcommand{\N}{\mathbb{N}}
\begin{document}

\title{Three geometric applications of quandle homology}

\author{Maciej Niebrzydowski}
\address{University of Louisiana at Lafayette,
Lafayette, LA 70504}
\email{mniebrz@gmail.com}

\date{July 4th, 2008}

\thispagestyle{empty}

\begin{abstract}

In this paper we describe three geometric applications of quandle homology. We show that it gives obstructions to tangle embeddings, provides the lower bound for the $4$-move distance between links, and can be used in determining periodicity of links. 
 
\end{abstract}

\maketitle

\section{Definitions and preliminary facts}  

\begin{definition}
A {\it quandle}, $X$, is a set with a binary operation 
$(a, b) \mapsto a * b$
such that
\begin{enumerate}
\item For any $a \in X$, 
$a* a =a$.

\item For any $a,b \in X$, there is a unique $c \in X$ such that 
$a= c*b$. 

\item For any $a,b,c \in X$,
$ (a*b)*c=(a*c)*(b*c)$ (right distributivity). 
\end{enumerate}
\end{definition}

Note that the second condition can be replaced with the following requirement:
the operation $*b\colon Q\to Q$, defined by $*b(x)=x*b$, is a bijection. The inverse map to $*b$ is denoted by
$\overline{*}b$. 

\begin{definition}
A {\it rack} is a set with a binary operation that satisfies 
conditions (2) and (3) from the definition of quandle.
\end{definition}

The following are some of the most commonly used examples of quandles.
\begin{itemize}
\item[-]
Any group $G$ with conjugation 
as the quandle operation:\\ $a*b=b^{-1} a b$. 
\item[-]
Let $n$ be a positive integer.
For elements  $i, j \in \{ 0, 1, \ldots , n-1 \}$, define
$i\ast j \equiv 2j-i \pmod{n}$.
Then $\ast$ defines a quandle
structure  called the {\it dihedral quandle},
  $R_n$.
It can be identified with  the
set of reflections of a regular $n$-gon
with conjugation
as the quandle operation.
\item[-]
Any $\mathbb{Z}[t, t^{-1}]$-module $M$
is a quandle with
$a*b=ta+(1-t)b$, for $a,b \in M$, called the {\it  Alexander  quandle}.
Moreover, if $n$ is a positive integer, then
$\mathbb{Z}_n[t, t^{-1}]/(h(t))$
is a quandle for
a Laurent polynomial $h(t)$.
\end{itemize}

The last example can be vastly generalized \cite{Joy}; for any group $G$ and its automorphism $\tau\colon G\to G$, $G$ becomes a quandle when equipped with the operation $g*h=\tau(gh^{-1})h$. If we consider the anti-automorphism $\tau(g)=g^{-1}$, we obtain another well known quandle, $Core(G)$, with
$g*h=hg^{-1}h$. 

Very likely the earliest work on racks is due to J. Conway and G. Wraith \cite{CW, FR}, who studied the conjugacy operation in a group.
The notion of quandle was introduced independently by D. Joyce \cite{Joy} and S. Matveev \cite{Mat}.

Joyce introduced the fundamental knot quandle, that is a classifying invariant of classical knots up to orientation-reversing
homeomorphism of topological pairs \cite{Joy}. However, just like in the case of fundamental groups, it is very hard to decide whether two given knot quandles are isomorphic. There are several other knot invariants derived from quandles that are easier to work with. For example, one can consider the family of all homomorphisms from the fundamental knot quandle to the given quandle, i.e., the set of all quandle colorings. The cardinality of this set is a knot invariant. 

\begin{figure}
\begin{center}
\includegraphics[height=3 cm]{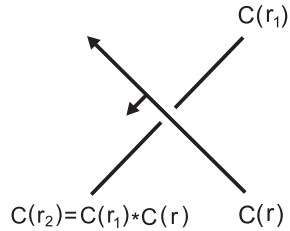}
\caption{The rule of quandle coloring at a crossing.}\label{crossingrule}
\end{center}
\end{figure}

\begin{definition}[\cite{CKS}]\label{colorings}
Let $X$ be a fixed quandle. Let $K$ be a given diagram of an oriented classical link, and let $R$ be the set of over-arcs of the diagram. The normals to arcs are given in such a way that the pair (tangent, normal) matches the usual orientation of the plane. A quandle coloring $C$ is a map $C\colon R\to X$ such that at every crossing, the relation depicted 
in Fig.\ref{crossingrule} holds. More specifically, let $r$ be the over-arc at a crossing, and $r_1,\ r_2$ be under-arcs such that the normal of the over-arc points from $r_1$ to $r_2$. Then it is required that 
$C(r_2)=C(r_1)*C(r)$.
\end{definition}

\begin{figure}
\begin{center}
\includegraphics[height=8 cm]{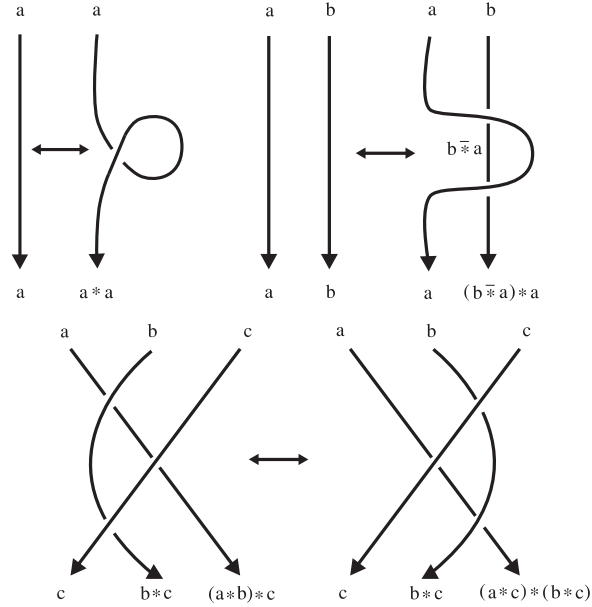}
\caption{Reidemeister moves and quandle axioms.\label{rad}}
\end{center}
\end{figure}

The axioms for a quandle correspond to the Reidemeister moves via quandle colorings of knot diagrams.
This correspondence is illustrated in Fig.\ref{rad}.

\section{Rack and quandle homology}\label{homologydef}

Rack homology and homotopy theory were first defined and 
studied in \cite{FRS}, and a modification to quandle homology theory 
was given in \cite{CJKLS} to define knot invariants in a state-sum form (so-called cocycle knot invariants).

Here we recall the definition of rack, degenerate and quandle homology
after \cite{CKS}.
\begin{definition} 
\begin{enumerate}
\item[(i)]
For a given rack $X$, let $C^R_n(X)$ be the free abelian 
group generated by $n$-tuples $(x_1,x_2,\ldots,x_n)$ of elements of $X$; 
in other words, $C^R_n(X) = {\Z}X^n = ({\Z}X)^{\otimes n}$.\\ 
Define a boundary   
homomorphism $\partial: C^R_n(X) \to C^R_{n-1}(X)$ by:
$$\partial(x_1,x_2,\ldots,x_n) = $$ 
$$\sum_{i=2}^n (-1)^i((x_1,\ldots,x_{i-1},x_{i+1},\ldots, 
x_n) - (x_1*x_i,x_2*x_i,\ldots,x_{i-1}*x_i,x_{i+1},\ldots,x_n)).$$ 
$(C^R_*(X),\partial)$ is called the rack chain complex of $X$.
\item[(ii)] Assume that $X$ is a quandle. Then there is a subchain 
complex $C^D_n(X) \subset C^R_n(X)$, generated by $n$-tuples $(x_1,\ldots,x_n)$ 
with $x_{i+1}=x_i$ for some $i$. The subchain complex $(C^D_n(X),\partial)$ 
is called the degenerated chain complex of a quandle $X$.
\item[(iii)] The quotient chain complex $C^Q_n(X)=C^R_n(X)/C^D_n(X)$ is 
called the quandle chain complex. 
\item[(iv)] The (co)homology  of rack, degenerate, and quandle chain complexes 
is called rack, degenerate, and quandle (co)homology, respectively.
\item[(v)] For an abelian group $G$, define the chain complex\\ $C^Q_*(X; G)=C^Q_*\otimes G$, with $\partial=\partial\otimes id$.
The groups of cycles and boundaries are denoted respectively by $ker(\partial)=Z^Q_n(X; G)\subset C^Q_n(X; G)$
and $Im(\partial)=B^Q_n(X; G)\subset C^Q_n(X; G)$.
The $n$th quandle homology group of a quandle $X$ with coefficient group $G$ is defined as
$$H^Q_n(X; G)=H_n(C^Q_*(X; G))=Z^Q_n(X; G)/B^Q_n(X; G).$$
\end{enumerate}
\end{definition}

Rack homology and quandle homology were studied by many authors, for example in \cite{CES, CJKLS, CJKS, CKS, EG,
FRS, LN, Moc}. Free part of rack (and quandle) homology is known for a large class of racks and quandles
\cite{EG, Moc}. However, there are many open problems concerning the torsion part. 

In this paper we will show how to use the information about homology of quandles in solving some geometric problems concerning knots and links. The effectiveness of these methods grows together with better understanding of quandle homology.

\section{Application to tangle embeddings}

First, we will explain, following \cite{Gr,CKS03,CKS01}, the procedure of assigning a cycle in quandle homology
to an oriented colored link diagram. 2-cycles correspond to diagrams with the usual quandle coloring, and
3-cycles are assigned to diagrams with shadow colorings.

\begin{figure}
\begin{center}
\includegraphics[height=4 cm]{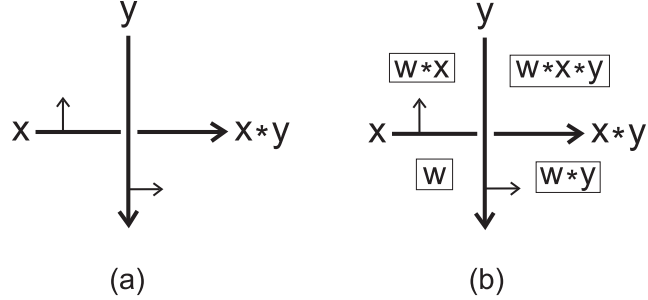}
\caption{Coloring and shadow coloring.\label{shadows}}
\end{center}
\end{figure}

\begin{definition}\label{shadowcol}
Let $Q$ be a fixed quandle, $D$ be a link diagram, and $\tilde{R}$ be the set of arcs and regions separated by the underlying immersed curve of $D$. A shadow coloring of $D$ is a function
$\tilde{C}\colon \tilde{R}\to Q$ satisfying the following two conditions.
\begin{itemize}
\item[(1)] The rules of labeling of arcs are as in the ordinary quandle coloring;
\item[(2)] Coloring of regions satisfies the condition illustrated in Figure \ref{shadows}(b), i.e.,
if $R_1$ and $R_2$ are two regions separated by an arc $r$ colored by $x$, and the normal vector to $r$ points from
$R_1$ to $R_2$, then the color of $R_2$ must be equal to $w*x$, where $w$ is the color of $R_1$.
\end{itemize}
Note that despite the fact that near the crossing there is more than one way to go from one region to another,
the third quandle axiom (the right distributivity) guarantees unique colors near a crossing. 
\end{definition}  

Let $D$ be a link diagram colored with elements of a finite quandle $X$.
Each positive crossing represents a pair $(x,y)\in C^Q_2(X)$, where $x$ is the color of an under-arc away 
from which points the normal of the over-arc labeled $y$ (see Figure \ref{shadows}(a)). 
In the case of negative crossing, we write $-(x,y)$.
The sum of such 2-chains taken over all crossings of the diagram forms a 2-cycle (see \cite{CKS03} for details).
Thus, it represents an element in $H^Q_2(X)$.

In the case of shadow coloring, each positive crossing corresponds to the triple 
\mbox{$(w,x,y)\in C^Q_3(X)$},
where $w$ is the color of so-called source region. It is the region near the crossing such that both normal vectors to the 
arcs colored by $x$ and $y$ point away from this region (see Figure \ref{shadows}(b), where colors assigned to the regions are depicted
as letters enclosed within squares).
A negative crossing represents the triple $-(w,x,y)$.
The sum of such signed triples taken over all crossings of $D$ gives an element of $H^Q_3(X)$ (\cite{CKS03}).
 
\begin{figure}
\begin{center}
\includegraphics[height=6.5 cm]{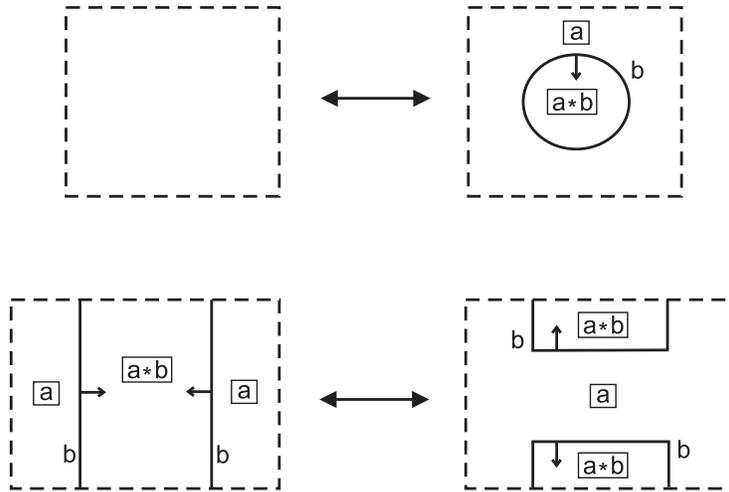}
\caption{Moves on shadow-colored diagrams that do not change homology class.\label{hommoves}}
\end{center}
\end{figure}   

Carter, Kamada, and Saito gave a list of moves on colored or shadow-colored link diagram that do not change the homology class 
represented by this diagram (\cite{CKS01, CKS}).
Their list includes Reidemeister moves and two moves illustrated in the Figure \ref{hommoves}. We are going to use these moves in our construction of obstructions to tangle embeddings.
The first move is creating or deleting a trivial component with appropriate shadow coloring.
The second move allows to change the connections between arcs of the diagram, if these arcs have the same color and 
opposite orientation.

A $2n$-tangle consists of $n$ disjoint arcs in the 3-ball. We ask the following question, 
that was first considered by D. Krebes \cite{Kre}. For a given knot $K$, and a tangle 
$T$, can we embed $T$ into $K$? In other words, is there a diagram of $T$ that extends to 
a diagram of $K$?
This problem is important due to its applications in the study of DNA. 
A number of knot invariants have been 
used to find criteria for tangle embeddings (see for example \cite{PSW}, \cite{Rub}).

Let us recall the definition of a special type of colorings of tangles that will be essential for defining homological obstructions to tangle embeddings.

\begin{definition}\label{bchrom}
Let $D_T$ be a tangle diagram, and $Q$ be a quandle. A 
\textit{boundary-monochromatic coloring} of $D_T$ is a map from the set of arcs of $D_T$ to quandle $Q$ satisfying the usual conditions for quandle colorings of knot diagrams, and an additional requirement that all boundary points receive the same color. 
\end{definition}

If a tangle $T$ embeds into a knot $K$, then each boundary-monochromatic coloring of $D_T$ can be extended trivially to the whole diagram of $K$, i.e., all arcs outside $D_T$ receive the color of the boundary points of $D_T$. Thus, the existence of nontrivial boundary-monochromatic colorings of $D_T$ gives the first basic obstruction to tangle embeddings, for $T$ can possibly embed only into knots admitting at least the same number of nontrivial colorings (see also \cite{Kre}).

\begin{definition}
A boundary-monochromatic shadow coloring of a tangle diagram $D$ is obtained from the ordinary boundary-monochromatic coloring of $D$ by choosing a color of any region of $D$ and extending this coloring to other regions according to the rules of Definition \ref{shadowcol}. Notice that such extension is unique.
\end{definition} 

\begin{lemma}
Every boundary-monochromatic coloring of an oriented diagram $D$ of a tangle $T$ with elements of a fixed quandle $X$ represents
an element in $H^Q_2(X)$. Every boundary-monochromatic shadow coloring of $D$ represents an element in $H^Q_3(X)$.
\end{lemma} 
\begin{proof}
The fact that all boundary points of $D$ have the same color allows us to take any closure of a diagram $D$ and obtain a colored link diagram that represents an element in $Z^Q_2(X)$ (or in $Z^Q_3(X)$ in the case of shadow coloring).
Any two such closures can be transformed one into another by a sequence of homology moves illustrated in the 
Figure \ref{hommoves}. Therefore, $D$ (as well as $T$) represents an element in quandle homology, i.e., element represented by any of its closures.
\end{proof}

Now we can define obstructions to tangle embeddings using quandle homology.

\begin{theorem}\label{homembeddings}
If a tangle $T$ embeds into a link $L$ then for every boundary monochromatic (shadow) coloring $\alpha$ of a diagram $D$ 
of $\,T$ there exists a 
(shadow) coloring $\beta$ of any diagram $\overline{D}$ of $L$, that represents the same homology class as the one represented by $\alpha$.
\end{theorem}
\begin{proof}
If $T$ embeds in $L$, then there exists a diagram $\widetilde{D}$ of $L$ such that $D$ is a part of it. Any boundary-monochromatic (shadow) coloring of $D$ extends trivially to a coloring
$\beta$ of $\widetilde{D}$.
Then, using homology moves from Figure \ref{hommoves}, one can destroy all crossings in $\widetilde{D}$ that are outside of 
$D$, and remove trivial components that may appear during this process. As a result one obtains one of the closures of $D$.
Homology class does not depend on the closure. Therefore, cycle represented by $\alpha$ equals to the cycle represented by 
$\beta$ in $H^Q_2(X)$ (or $H^Q_3(X)$ in the case of shadow colorings). Finally, any coloring of $\widetilde{D}$ gives a coloring of any other diagram $\overline{D}$ of $L$ by a sequence of Reidemeister moves (they do not change the homology class). 
\end{proof}

\begin{figure}
\begin{center}
\includegraphics[width=11 cm]{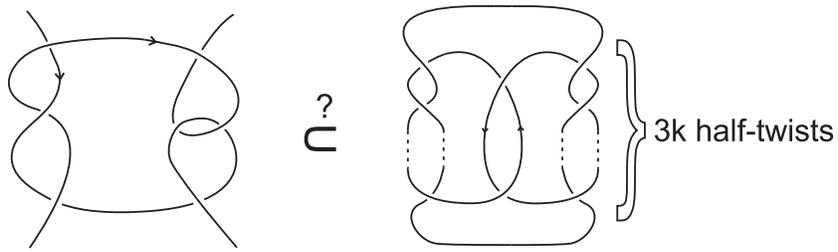}
\caption{An example of a problem of embedding a given tangle into a link.\label{tangleinknot}}
\end{center}
\end{figure}  

\begin{figure}
\begin{center}
\includegraphics[height=9.5 cm]{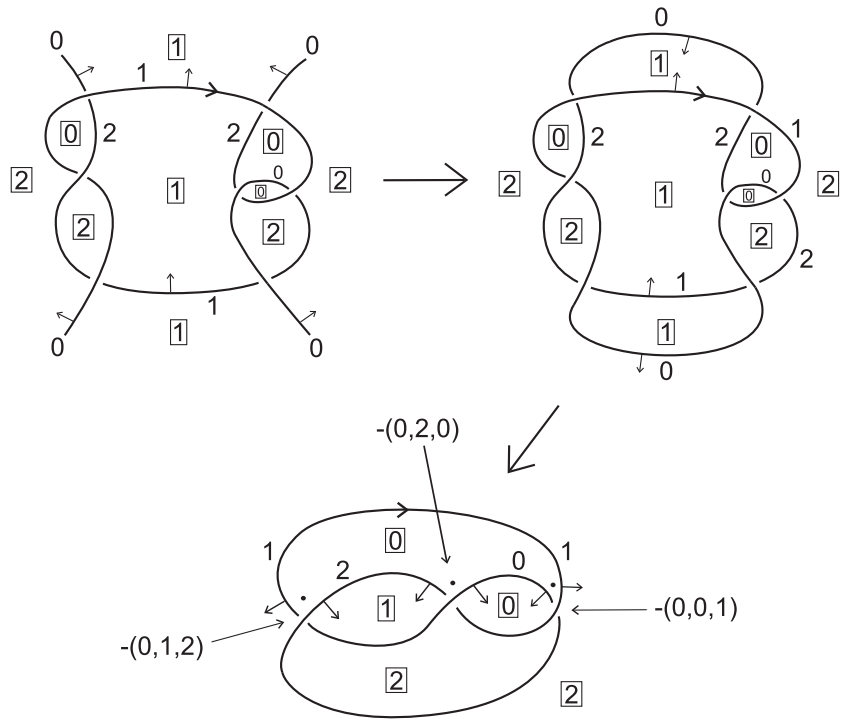}
\caption{A shadow-colored tangle that represents a generator of $H^Q_3(R_3)$.\label{tanglenot}}
\end{center}
\end{figure}

\begin{figure}
\begin{center}
\includegraphics[height=7 cm]{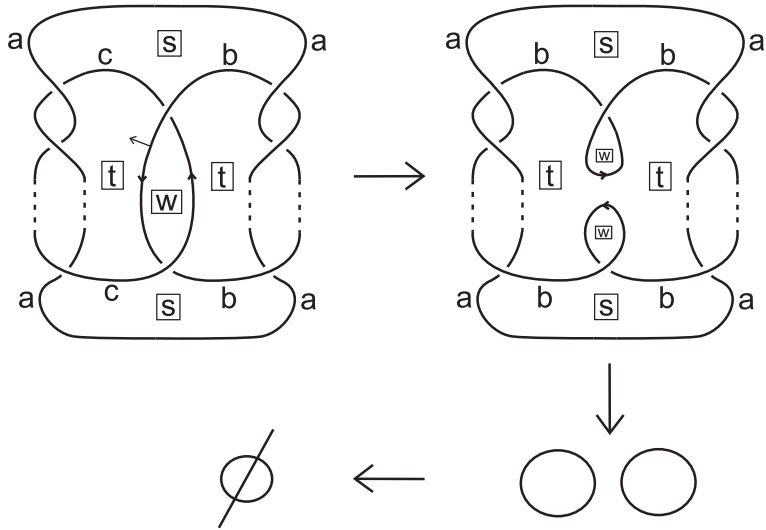}
\caption{A nontrivial knot, whose every coloring represents 0 in $H^Q_3(R_3)$.\label{knotw}}
\end{center}
\end{figure}  

\noindent\textbf{Example.} Figure \ref{tangleinknot} illustrates an example of a tangle $T$, and
a family of links that have $3k, k\in\N,$ half-twists on each side. Let $L$ denote any member of this family.
We can use the third quandle homology of the dihedral quandle $R_3$ to show that $T$ does not embed into $L$.
Figure \ref{tanglenot} shows an example of a boundary-monochromatic shadow coloring of a diagram of $T$ with elements of $R_3$.
The nominator closure of this tangle is a shadow-colored trefoil knot. This coloring
represents a chain $-(0,1,2)-(0,0,1)-(0,2,0)$ that gives a generator of $H^Q_3(R_3)$ \cite{NP2}.
On the other hand, every coloring of the link $L$ represents $0$ in $H^Q_3(R_3)$. We can see it as follows.
Quandle $R_3$ is the simplest nontrivial example of a Burnside kei \cite{NP1}, and is invariant under 3-moves and, more generally,
under $3k$-moves. That is why we can label the top arcs and the corresponding bottom arcs of the diagram of $L$ with the same elements $a$, $b$, $c$.
Such labeling forces relations $b*c=b$ and $c*b=c$ that imply the equality $b=c$ in $R_3$.
Because of this equality it is possible to perform a homology move on the diagram of $L$ (see Figure \ref{tangleinknot}),
that transforms it into unlink representing $0$ in homology.
From the Theorem \ref{homembeddings} follows that $T$ cannot be embedded into $L$. 

Another quandle-based approach to the tangle embedding problem, using quandle cocycle invariants, was proposed in \cite{AERSS}. 

\section{The structure of $H^Q_2(R_4)$}\label{struc}
In order to provide examples for the next two applications, we will now analyze the second homology group of the dihedral quandle $R_4$. 

Let us recall that the dihedral quandle $R_4$ is a set $\{0, 1, 2, 3 \}$ with operation
$i * j \equiv 2j-i \pmod{4}$. It consists of two orbits (with respect to the action of $R_4$ on itself by the right multiplication): $\{0,2\}$ and $\{1, 3\}$. 

To simplify our notation, and make it more general, we will write the elements of this quandle as
$\{a, b, a*b, b*a\},$ where $a$ and $b$ are representatives of different orbits.
Note that the elements of $R_4$ (when written as longer products involving $a$ and $b$) can be determined by looking at the first letter in the word, and the parity of the letter from the second orbit that appears in the rest of the word. For example, if $s$ denotes the  number of $b$'s in the word, then

\begin{displaymath}
a*b*\ldots=
\left\{ \begin{array}{l}
a,\ \textrm{if $s$ is even}\\
\ \ \ \ \ \ \ \\
a*b,\ \textrm{if $s$ is odd}\\  
\end{array} \right.  
\end{displaymath}

It is known (see for example \cite{LN}) that $H^Q_2(R_4)=Z^2\oplus (Z_2)^2.$
We will show that the free part is generated by:
\begin{enumerate}
\item[] $f_1=(a,b)+(a*b,b),$ 
\item[] $f_2=(b,a)+(b*a,a),$
\end{enumerate}
and that the generators of the torsion part are:
\begin{enumerate}
\item[] $t_1=(a,a*b),$ 
\item[] $t_2=(b,b*a).$
\end{enumerate}
The first part of the statement follows from evaluating the cocycles 
\begin{enumerate}
\item[] $\chi_{(a,b)}+\chi_{(a,b*a)},$ 
\item[] $\chi_{(b,a)}+\chi_{(b,a*b)},$
\end{enumerate}
on $f_1$ and $f_2$. Here, $\chi_{(x,y)}$ denotes the characteristic function of $(x,y)$, and the above cocycles were proven to be generators of 
$H^2_Q(R_4,\Z)=\Z^2$ in \cite{CJKLS}.
To prove the second part, we first notice that $t_1$ and $t_2$ are either torsion elements or $0$, since we have
$$\partial(-(a*b,b,a*b)-(a*b,b*a,a*b))=2 t_1,$$
$$\partial(-(b*a,a,b*a)-(b*a,a*b,b*a))=2 t_2.$$
To prove non-triviality, we use the following cocycles $c_1$ and $c_2\in H^2_Q(R_4,\Z_2)$:
\begin{enumerate}
\item[] $c_1=\chi_{(a,b)}+\chi_{(a*b,b)}+\chi_{(a,a*b)}+\chi_{(a*b,a)},$ 
\item[] $c_2=\chi_{(b,a)}+\chi_{(b*a,a)}+\chi_{(b,b*a)}+\chi_{(b*a,b)}.$
\end{enumerate}
Since $c_1(t_1)=1$ and $c_2(t_2)=1$, $t_1,t_2$ and $c_1,c_2$ must be non-trivial. We also note that $c_1$ and $c_2$ evaluate trivially on $f_1$ and $f_2$.

\section{The lower bound for the $4$-move distance between links}

\begin{figure}
\begin{center}
\includegraphics[height=2 cm]{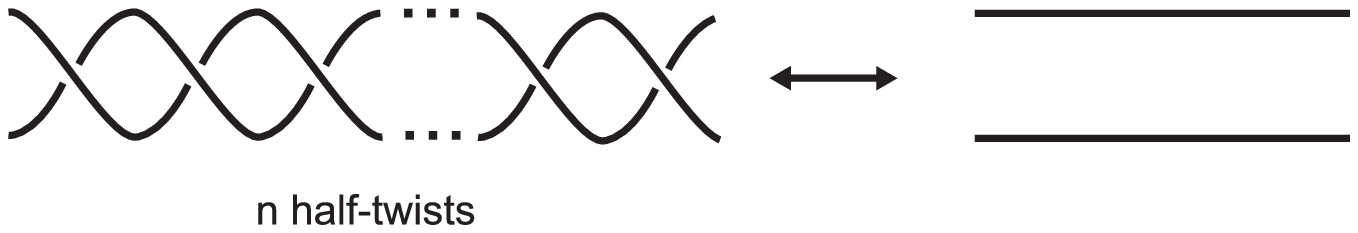}
\caption{$n$-move.\label{nmove}}
\end{center}
\end{figure}

\begin{definition}
An $n$-move is a replacement of $n$ half-twists by two parallel strings or vice versa in a link diagram (see Fig.\ref{nmove}).
\end{definition}

Of particular interest in knot theory are $4$-moves. One of the reasons is the following old conjecture \cite{Kir,Pr1,Pr2}.

\begin{conjecture}[Nakanishi, 1979]
Every knot is $4$-move equivalent to the trivial knot. In other words, every knot can be transformed into a trivial knot using $4$-moves and Reidemeister moves.
\end{conjecture}

Not every link is $4$-move equivalent to a trivial link, in particular, the linking
matrix modulo $2$ is preserved by $4$-moves. Furthermore, Nakanishi demonstrated
that the Borromean rings cannot be reduced to the trivial link of three components \cite{Nak,Pr2}.
Kawauchi expressed the question for links as follows:
\begin{problem}[\cite{Kir}]
\begin{itemize}
\item[]
\item[(i)] Is it true that if two links are link-homotopic then they are $4$-move equivalent?
\item[(ii)] In particular, is it true that every $2$-component link is $4$-move equivalent
to the trivial link of two components or to the Hopf link?
\end{itemize}
\end{problem}

A $3$-component counterexample to this problem was provided in \cite{DP}. The second part of the question remains open, and is actively investigated. 

In this paper we consider the following problem.

\begin{problem}
If two links are $4$-move equivalent, what is the minimal number of $4$-moves needed to transform one into the other?
\end{problem}

Therefore, it is natural to make the following definition.

\begin{definition}
An $n$-move distance, $d_n(L_1,L_2)$, between two links 
$L_1$ and $L_2$, is the minimal number of $n$-moves realizing the $n$-move equivalence, or $\infty$ if $L_1$ and $L_2$ are not $n$-move equivalent.
\end{definition}

For example, the $4$-move distance between the trivial link of two components and the Hopf link is $\infty$, as indicated by their linking matrices modulo $2$.

\begin{figure}
\begin{center}
\includegraphics[height=2 cm]{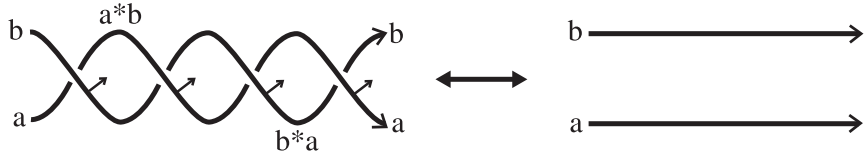}
\caption{A colored $4$-move represents a cycle in $H^Q_2(R_4)$.\label{czmove}}
\end{center}
\end{figure}

We will now explain how the quandle homology of $R_4$ can be often used to obtain the lower bound for $d_4(L_1,L_2)$.

\begin{lemma}\label{twistcycle}
Any $R_{2k}$-coloring of the two oriented strings with $2k$ half-twists represents a cycle in $H^Q_2(R_{2k})$.
\end{lemma}
\begin{proof}
It is known (see for example \cite{Pr1,Pr2,NP1}) that the dihedral quandle $R_n$ is an invariant under $n$-moves. In other words, for any 
$R_n$-coloring of the two strings with $n$ half-twists, the colors of the two initial arcs are the same as colors of the corresponding final arcs (see Figure \ref{czmove} for an illustration of this fact in the case of $4$-move). Colorings with dihedral quandles do not depend on the orientation of the link. However, if we want to analyze quandle homology, the orientation has to be taken into account. In the case of an even number of half-twists, for any orientation (parallel or anti-parallel) of the twisted strings, it is possible to join the upper left arc with the upper right arc, and the lower left arc with the lower right arc, without introducing any additional crossings. Thus, we obtain a properly colored and oriented, uniquely determined link that represents an element in $H^Q_2(R_{n})$. That is not the case when $n$ is odd and the orientation is anti-parallel. Often the chain determined by such $n$ colored crossings is not a cycle. 
\end{proof}

\begin{lemma}\label{moves}
Let $c$ be a cycle representing some $R_4$-coloring of the two oriented strings with $4$ half-twists. If both strings have colors from a single orbit, then $c$ is homologically trivial, otherwise $c=\pm (f_1+f_2+t_1+t_2)$, where $f_1$, $f_2$, $t_1$, $t_2$ are as in the previous section.
\end{lemma}

\begin{proof}
First, for a given quandle $X$, we define the map
$$*_a: C^Q_{n}(X) \to C^Q_{n}(X)$$ determined by $*_a(w)=w*a$, for any $w\in X^n$, or more precisely, 
$$*_a(x_1,\ldots,x_n) = (x_1,\ldots,x_n)*a= (x_1*a,\ldots,x_n*a).$$
We are going to use the following fact from \cite{NP2}: if $z$ is a cycle, then $*_a(z)$ is also a cycle, homologous to $z$, for any $a\in X$.

Applying this fact to $t_1$ and $t_2$, we can check that any pair of different elements of $R_4$ from the same orbit represents a torsion in $H^Q_2(R_{4})$.
It follows that if the strings are nontrivially colored by elements from the same orbit, then such coloring represents a cycle homologous to $\pm 4t_1$ or $\pm 4t_2$. It can be checked by inspection that each coloring that uses the elements from different orbit gives a cycle $c$ that decomposes into two smaller cycles: $c=\pm (c_1+c_2)$, where $c_1=(a,b)+(a*b,b*a)$ and $c_2=(b,a)+(b*a,a*b)$. Sometimes cycles $c_1'=(a*b,b)+(a,b*a)$ or $c_2'=(b*a,a)+(b,a*b)$ appear (as in the Fig.\ref{czmove}), but $c_1=*_b(c_1')$ and $c_2=*_a(c_2')$, so there is no difference in homology. Finally, we check that
$$\partial((a,b,a*b)-(a,b,b*a))=c_1-(f_1-t_1)$$
$$\partial((b,a,b*a)-(b,a,a*b))=c_2-(f_2-t_2),$$
and, since homologically $t_1=-t_1$ and $t_2=-t_2$, the proof is finished. 
\end{proof}

\begin{figure}
\begin{center}
\includegraphics[height=6.5 cm]{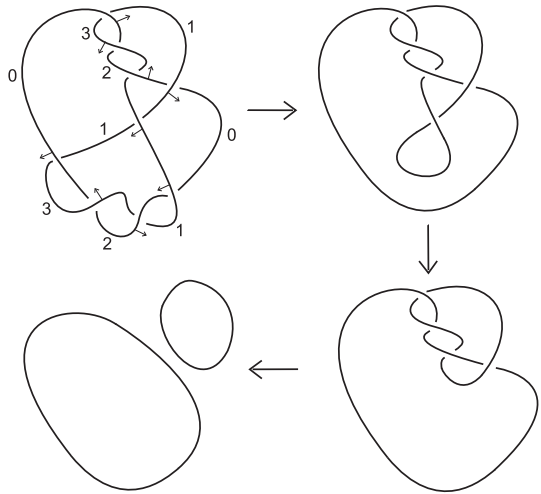}
\caption{A link reduced to the trivial link by two $4$-moves.\label{moves4ex}}
\end{center}
\end{figure}

\begin{corollary}
Let $L$ be an oriented link diagram colored with elements of the dihedral quandle $R_4$, and let $c$ be a cycle in $H^Q_2(R_{4})$ represented by this coloring. Then, for any $4$-move performed on the diagram, $c$ either remains unchanged or is replaced by $c\pm (f_1+f_2+t_1+t_2)$.
\end{corollary}

\begin{corollary}
The $4$-move distance, $d_4(L_1,L_2)$, between links $L_1$ and $L_2$, such that at least one of them admits nontrivial $R_4$-colorings, can be analyzed by comparing the multiplicities of $f_1+f_2$ appearing in the cycles represented by the colorings.
\end{corollary}

\noindent\textbf{Example.} We consider the coloring of the link illustrated in Figure \ref{moves4ex}. It represents a cycle $c\in H^Q_2(R_{4})$ of the form:
$$(a,b)+(a*b,b*a)+(b,a)+(b*a,a*b)+(a,b*a)+(a*b,b)+(b*a,a)+(b,a*b).$$ Using a similar technique as in the proof of Lemma \ref{moves}, we can conclude that 
$$c=2(f_1+f_2+t_1+t_2)=2(f_1+f_2).$$ Thus, at least two $4$-moves are necessary to reduce it to the trivial link with two components (whose colorings represent $0$ in homology). It cannot be reduced to the Hopf link, since the Hopf link admits only colorings using elements of one orbit, and this property is preserved by $4$-moves. As shown in the Figure \ref{moves4ex}, two $4$-moves suffice to make the reduction.

We note that any link with a coloring representing a cycle that is not a multiple of $f_1+f_2+t_1+t_2$ would be a counterexample to the second part of   
Kawauchi's question, because the colorings of the Hopf link and the colorings of the trivial link do not give any nontrivial classes in $H^Q_2(R_{4})$. 
No such link has been found so far. However, since every cycle from the second homology can be represented by a colored virtual link (\cite{CKS}), above technique provides virtual counterexamples to the question. One such example is a virtual link with a coloring representing $f_1+f_2+t_1$.

We also remark that the above method can be generalized to quandles $R_{2k}$ and the $2k$-move distance. More generally, it should work with certain rational moves (see \cite{DP} for a definition) and the rational move distance.

\section{Application to periodicity of links}

\begin{figure}
\begin{center}
\includegraphics[height=4.5 cm]{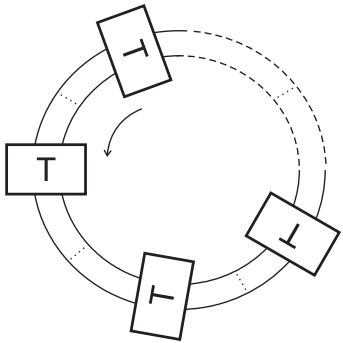}
\caption{Invariant diagram of a periodic link.\label{invdiag}}
\end{center}
\end{figure}

\begin{definition}
Let $p$ be a prime number.
A link $L$ in $S^3$ is called $p$-periodic if there is a $Z_p$-action on $S^3$, with a circle as a fixed point set, which maps $L$ onto itself,
and such that $L$ is disjoint from the fixed point set. Furthermore, if $L$ is an
oriented link, one assumes that each generator of $Z_p$ preserves the orientation
of $L$ or changes it to the opposite one.
\end{definition}
By the positive solution to the Smith Conjecture (\cite{MB}), if a link is $p$-periodic, then $L$ has a diagram 
$\widetilde{D}$ such that the rotation by an angle $\frac{2\pi}{p}$ about a point away from the diagram leaves
$\widetilde{D}$ invariant. There exists an $n$-tangle $T$ such that $L$ is the closure of $T^p$, i.e., a tangle obtained by gluing $p$ copies of $T$ in a natural way, as illustrated in Figure \ref{invdiag} (see also \cite{CL,GKP,PS}).

In this section we will show that sometimes we can use quandle homology to prove that a link $L$ is
not $p$-periodic for some prime $p$.

\begin{theorem}
Let $p$ be a prime number, and $D$ be a diagram of a $p$-periodic link $L$.
If a coloring (shadow coloring) of $D$ with elements of some fixed quandle $X$ represents a homology class $c$ in $H^Q_2(X)$
(or $H^Q_3(X)$ in the case of shadow coloring), then either there exist $p-1$ different colorings of $D$
that represent the same element in homology as $c$, or $c=p\,\widetilde{c}\,$, for some $\widetilde{c}$ in $H^Q_2(X)$ 
(or in $H^Q_3(X)$).
\end{theorem}

\begin{figure}
\begin{center}
\includegraphics[height=5 cm]{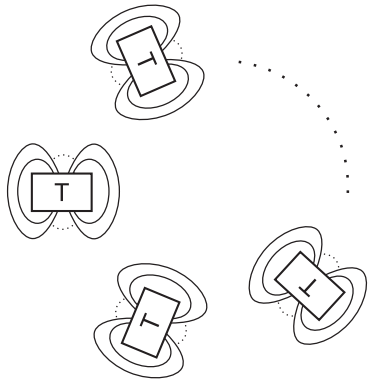}
\caption{Homological decomposition of a colored invariant diagram.\label{homdecomp}}
\end{center}
\end{figure}

\begin{proof}
First, let us note that there are two types of colorings of the invariant diagram $\widetilde{D}$ of a $p$-periodic
link $L$. One possibility is that all tangles $T$ that are building blocks of $\widetilde{D}$ receive exactly the same coloring. Otherwise, there is some asymmetry in the coloring, and it is possible to obtain from it $p-1$ different colorings 
by rotating the given coloring by a multiple of an angle $\frac{2\pi}{p}$. Note that in this process the position of the link diagram $\widetilde{D}$ is not changed, only coloring is rotated. This distinction becomes more clear if we translate each such coloring (via Reidemeister moves) into a coloring of some other, less symmetric diagram of $L$.
If $p$ is not prime, then we might obtain a smaller number of colorings than $p-1$, because some of them may be identical.
Let $C$ be any (shadow) coloring of a diagram $D$. If $C$ is of the first type, then using homology move that changes
connections between strings with the same color (see Figure \ref{hommoves}), we can decompose colored diagram $D$ into $p$ identical smaller diagrams $\widehat{D}$ (as in Figure \ref{homdecomp}). In this case, the element in quandle homology that is represented
by the coloring $C$ is equal to $p\,\widehat{C}$, where $\widehat{C}$ is element of homology corresponding
to $\widehat{D}$. If the coloring is of the second type, then each of the $p-1$ colorings obtained by rotating the original
coloring represents the same element in homology.
\end{proof} 

\begin{figure}
\begin{center}
\includegraphics[height=6 cm]{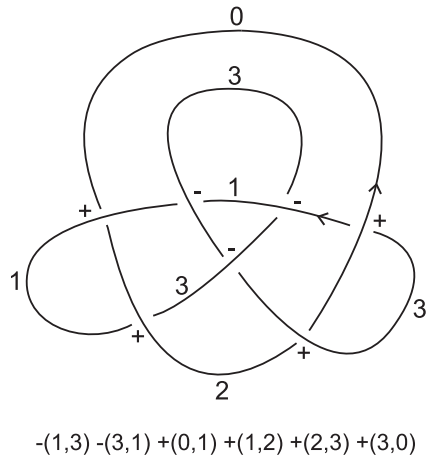}
\caption{A quandle coloring of the link $7^2_5$ and the nontrivial cycle it represents.\label{link}}
\end{center}
\end{figure} 

\noindent\textbf{Example.}
We can use the second quandle homology of the quandle $R_4$ to show that $2$ is the only possible period for the link $7^2_5$.
This link has $16$ colorings using the quandle $R_4$. Eight of them are either trivial or represent $4t$, where $t$ is an element from the $\Z_2$-torsion.
The remaining eight colorings give cycles of the form:
$$(a,b)+(a*b,b*a)+(b,a*b)+(b*a,a)-(b,b*a)-(b*a,b).$$ As in the previous section, we can recognize them as homologous to $f_1+f_2+t_1+t_2$.
 The possibility of such element being equal to $p$ times some other element, for $p$ prime, is excluded. If the link $7^2_5$ were $p$-periodic, then the aforementioned $8$ colorings would have to be partitioned into $p$-element subsets. Therefore, the only candidate for the prime period of the link $7^2_5$ is 2.

\end{document}